# NONEXISTENCE OF COUNTINUOUS PEAKING FUNCTIONS

Jiye Yu

October 1995

ABSTRACT. We construct a smoothly bounded pseudoconvex domain such that every boundary point has a p.s.h. peak function but some boundary point admits no (local) holomorphic peak function.

## 1. Introduction and results

Let $\Omega$ be a domain in $\mathbb{C}^n$. Denote by $A(\Omega)$ the Banach algebra of continuous functions on $\overline{\Omega}$ which are holomorphic in $\Omega$. We say that a boundary point $p$ is a (continuous) *peak point* for $\Omega$ if there is a function $f \in A(\Omega)$ such that: $f(p) = 1$ and $|f(z)| < 1$ on $\overline{\Omega} \setminus \{p\}$. The function $f$ is called a *peaking function.* $p$ is said to be a *local peak point* if there is an open neighborhood $U$ of $p$ such that $p$ is a peak point for $\Omega \cap U$. The existence/nonexistence of peaking functions has many important consequences in the function theory of several complex variables, we refer the reader to [BF1-2] [B-P] [BSY] [CHR] [FOR] [KR2] [RA1] [SI1], to name only a few.

A basic problem associated with peak points is to determine a necessary and sufficient condition for a boundary point of a pseudoconvex domain to be a peak point. There are several known sufficient conditions [BF1] [F-S] [F-M] [NO1] [D-H] [YU]. For example, if $p \in b\Omega$ is strongly pseudoconvex point [RO1] or $p$ is *h-extendible* [YU], then it is peak point. On the other hand, the maximum principle implies that if $p \in b\Omega$ is a peak point, then the boundary of $\Omega$ cannot contain any nontrivial complex variety containing $p$ as an interior point. However, it is still unknown what kind of non-existence of analytic structures is necessary. Sibony constructed an example in [SI1] for which some boundary point is a local peak point but not a global one. Nevertheless, the boundary of his example contains complex discs. In this note, we would like to provide an example whose boundary contains no complex varieties, but not every boundary point is a peak point, not even a local peak point. As a matter of fact, our example domain will be B-regular which implies in a strong sense that there is no analytic structure in the boundary [SI2]. As far as smooth peak functions are concerned, there are already many counterexamples, see [FST] [K-N] [NO2].

1991 *Mathematics Subject Classification.* 32F15, 32F25.
Supported in part by NSF grant DMS-9500916, and at MSRI by NSF grant DMS-9022140.

Typeset by $\mathcal{AMS}$-TEX

Recall that a domain $\Omega$ is *B-regular* [SI2-3] [CA2] at $p \in b\Omega$ if there is a plurisubharmonic (**psh**) function $u$ on $\Omega$ continuous on $\overline{\Omega}$ such that $u(p) = 0$ and $u < 0$ on $\overline{\Omega} \setminus \{p\}$. Clearly, a peak point is necessarily B-regular. Moreover, the maximum principle excludes the existence of any analytic structure in the boundary around a B-regular point. B-regular domains are an important class of weakly pseudoconvex domains on which the $\overline{\partial}$-Neumann operator is well behaved (satisfies the compactness estimates) [CA2].

Our main results are

**Theorem 1.** *There is a smoothly bounded B-regular pseudoconvex domain in $\mathbb{C}^3$ such that some boundary point is not a peak point.*

As a consequence of the proof, we will also show that

**Theorem 2.** *There is a smoothly bounded B-regular pseudoconvex domain in $\mathbb{C}^3$ such that some boundary point is not a local peak point.*

*Acknowledgements:* The author is indebted to Professor Edgar Lee Stout for many stimulating conversations and helpful suggestions. He would also like to thank Professors Steven Krantz, John Wermer and Emil Straube for helpful discussions. This work was done when the author was visiting MSRI in Fall, 1995. The auther wishes to thank the institute for the support and hospitality.

## 2. The construction of the example

First of all, we fix some notations. Let $K$ be a compact subset of $\mathbb{C}^n$. As usual, $C(K)$ is the Banach algebra of countinuous functions on $K$. $\mathcal{O}(K)$ is the (uniform) closure in $C(K)$ of the functions holomorphic in some neighborhood of $K$. $PSH(K)$ is the closure in $C(X)$ of continuous psh functions in some neighborhood of $K$. $\mathcal{P}(K)$ is the closure in $C(K)$ of the set of polynomials in $\mathbb{C}^n$, and $\mathcal{R}(K)$ stands for the closure in $C(K)$ of the set of rational functions which have no poles in $K$.

**Lemma 1.** *There are a polynomially convex compact subset $K$ in $\mathbb{C}^2$ and $z_0 \in K$ satisfying*
  (i) *$K$ is B-regular, i.e. every point is a peak point for $PSH(K)$;*
  (ii) *There is a nontrivial representing measure for $z_0$ on $K$ for $\mathcal{O}(K)$.*

Here and below, "nontrivial" means "not a point mass".

In order to construct such a compact set, we first choose a compact set $X$ in $\mathbb{C}$ so that $\mathcal{R}(X)$ is *normal* but $\mathcal{R}(X) \neq C(X)$. Here "normal" means that for every pair of disjoint closed subsets $E_0, E_1$ of $X$, there is a rational function $r \in \mathcal{R}(X)$ such that $r \equiv i$ on $E_i$, $i = 0, 1$. Such a compact planar set was first found by McKissick [MCK] using a special "Swiss cheese" of the form $\overline{\Delta} \setminus \cup \Delta_j$, where each $\Delta_j$ is a small open disc inside the unit disc $\Delta$. In particular, $X \subset \overline{\Delta}$. By a theorem of Rossi [RO1] (see also [STO] [WER]), $\mathcal{R}(X)$ is doubly generated Banach algebra. More precisely, there are $r_1, r_2 \in \mathcal{R}(X)$, such that $\mathcal{R}(X) = \mathbb{C}[r_1, r_2]$. In fact, one can even choose $r_1(\zeta) = \zeta$.

Now define a map $\Phi : \mathbb{C} \to \mathbb{C}^2$ by $\Phi(\zeta) = (r_1(\zeta), r_2(\zeta))$. Then $\Phi$ is a holomorphic embedding of $X$ in $\mathbb{C}^2$. Denote its image by $K = \Phi(X)$. We will show that such $K$ satisfies the conditions in Lemma 1. Observe that the compact set $K$ is polynomially convex, i.e., for any $z \notin K$, there is a polynomial $P$ such that $|P(z)| > \|P\|_{C(K)}$. This follows from the fact that the joint spectrum of a set of generators for a commutative Banach algebra with identity is always polynomially convex and the fact that the spectrum of $\mathcal{R}(X)$ is $X$ itself (see, e.g. [STO] for details).

We first verify that $K$ is B-regular. Suppose not. Then in light of an equivalent condition for B-regularity given in [SI2], there are a point $z_0 \in K$ and a nontrivial Jensen measure $\mu$ for $z_0$ on $K$ for $PSH(K)$. That is, the measure $\mu$ is a probability measure on $K$ that satisfies

$$(1) \qquad u(z_0) \leq \int_K u(z)\, d\mu, \quad \forall u \in PSH(K).$$

Since $\mu$ is not a point mass, there must be two disjoint closed subsets $F_0, F_1$ of $K$ such that $\mu(F_i) > 0$, $i = 0, 1$. Set $E_i = \Phi^{-1}(E_i) \cap X$. Then $E_0, E_1$ are disjoint and compact. It follows from the normality of $\mathcal{R}(X)$ that there exists $r \in \mathcal{R}(X)$ such that $r \equiv i$ on $E_i$ for $i = 0, 1$. On the other hand, since $\mathcal{R}(X)$ is generated as a Banach algebra by $r_1, r_2$, there exist polynomials $g_n$ in $\mathbb{C}^2$ such that $r(\zeta) = \lim_{n \to \infty} g_n(r_1(\zeta), r_2(\zeta))$, uniformly on $X$. For any integers $m, n \geq 1$, the function $u_{mn} := \log(|g_n| + 1/m)$ is continuous and psh in $\mathbb{C}^2$. Applying (1) to $u_{mn}$, we obtain

$$\log(|g_n(z_0)| + \frac{1}{m}) \leq \int_K \log(|g(z)| + \frac{1}{m})\, d\mu.$$

If we set $\sigma = \mu \circ \Phi$, then $\sigma$ is a probability measure on $X$ such that $\sigma(E_i) = \mu(F_i) > 0$ for $i = 0, 1$. Let $\zeta_0 = \Phi^{-1}(z_0)$. A simple change of variables (see e.g. [HAL]) yields

$$\log(|g_n(\Phi(\zeta_0))| + \frac{1}{m}) \leq \int_X \log(|g_n(\Phi(\zeta))| + \frac{1}{m})\, d\sigma.$$

Letting $n \to \infty$ in the above inequality, we obtain

$$\log(|r(\zeta_0)| + \frac{1}{m}) \leq \int_X \log(|r(\zeta)| + \frac{1}{m})\, d\sigma.$$

Observe that the integrand in the above integral is bounded from above by the constant $\log(\|r\|_X + 1)$ and is non-increasing with respect to $m$. The monotone convergence theorem (see [HAL]) implies that

$$\log |r(\zeta_0)| \leq \int_X \log |r(\zeta)|\, d\sigma.$$

The right hand side of the above inequality is $-\infty$ due to the fact that $r(E_0) = \{0\}$ and $\sigma(E_0) > 0$. Thus we have $r(\zeta_0) = 0$. On the other hand, if we replace $g_n$ by $1 - g_n$ in the above arguments, we will also have $1 - r(\zeta_0) = 0$. This is clearly impossible. Consequently, $K$ must be B-regular.

In order to verify (ii) in Lemma 1 for $K$, we recall that for a compact planar set $Y$, $\mathcal{R}(Y) = C(Y)$ if and only if almost all points (in the Lebesgue measure) of $Y$ are peak points for $\mathcal{R}(Y)$ (see [STO]). Since for our $X$, $\mathcal{R}(X) \neq C(X)$, it follows that there exist at least one point $\zeta_0 \in X$ that is not a peak point for $\mathcal{R}(X)$. Necessarilly $\zeta_0 \in \Delta$. Thus by a theorem of Bishop [BIS] (see also [GAM]), there is a nontrivial representing measure $\sigma$ for $\zeta_0$ on $X$ for $\mathcal{R}(X)$. That is,

$$g(\zeta_0) = \int_X g(\zeta)\,d\sigma, \quad \forall g \in \mathcal{R}(X). \tag{3}$$

Set $\mu = \sigma \circ \Phi^{-1}$. Then $\mu$ is a well defined nontrivial probability measure on $K$. For any $f \in \mathcal{P}(K)$, let $g = f \circ \Phi \in \mathcal{R}(X)$. A simple application of change of coordinates for measurable mapping [HAL] to (3) yields

$$f(\Phi(\zeta_0)) = \int_X f(\Phi(\zeta))\,d\sigma = \int_K f\,d\mu.$$

It follows that

$$f(z_0) = \int_K f\,d\mu, \quad \forall f \in \mathcal{P}(K). \tag{4}$$

Here $z_0 = \Phi(\zeta_0)$. It remains to observe that $\mathcal{P}(K) = \mathcal{O}(K)$ by Oka-Weil theorem (see e.g. [RA2] [KR1]). Consequently, $\mu$ satisfies (ii) in Lemma 1. This finishes the proof of Lemma 1.

Next, we construct our domains in $\mathbb{C}^3$ to meet the requirements of Theorem 1 and 2. The construction follows from a theorem of Catlin [CA2]:

**Proposition 2.** *Let $K$ be any compact polynomially convex set in $\mathbb{C}^n$. Then there exists a $C^\infty$ smooth psh exhaustion function $\phi$ in $\mathbb{C}^n$ such that $\phi \geq 0$, $\phi^{-1}(0) = K$ and $\phi$ is strictly psh exactly in $\mathbb{C}^n \setminus K$.*

Let $\phi$ be the function obtained by applying Proposition 2 to our compact set $K$ in $\mathbb{C}^2$. Define a family of domains $\Omega_\delta$ for $\delta > 0$ in $\mathbb{C}^3$ by

$$\Omega_\delta = \{(z, w) \in \mathbb{C}^2 \times \mathbb{C} : |w|^2 + \phi(z) < \delta^2\}.$$

By Sard's theorem, for almost all $\delta > 0$, the domain $\Omega_\delta$ has $C^\infty$ boundary. Choose such a $\delta$. By replacing $\phi$ by $\phi/\delta^2$ and $w$ by $w/\delta$, we may assume that $\delta = 1$ and then set $\Omega = \Omega_1$ and $\mathbb{T} = \{\zeta \in \mathbb{C} : |\zeta| = 1\}$.

**Lemma 3.** *The domain $\Omega$ thus constructed is smoothly bounded pseudoconvex B-regular domain in $\mathbb{C}^3$. Moreover, $A(\overline{\Omega})|_{K \times \mathbb{T}} \subset \mathcal{O}(K \times \mathbb{T})$.*

Indeed, by Proposition 2, the set of weakly pseudoconvex boundary points of $\Omega$ is exactly $K \times \mathbb{T}$ which is B-regular by Lemma 1. It then follows that $\Omega$ is B-regular [SI2]. Next, for any $f \in A(\overline{\Omega})$, define a sequence of functions $f_n$ by $f_n(z, w) = f(z, nw/(n+1))$, $n \geq 1$. Then $f_n \to f$ uniformly on $\overline{\Omega}$ and $f_n$ is holomorphic on the domain

$$\{(z, w) \in \mathbb{C}^2 \times \mathbb{C} : \frac{n^2}{(n+1)^2}|w|^2 + \phi(z) < 1\} \supset \Omega,$$

which contains some open neighborhood of $K \times \mathbb{T}$. Thus Lemma 3 follows.

Now we prove Theorem 1. Let $p_0 = (z_0, 1) \in b\Omega$, where $z_0$ is given in Lemma 1. We claim that $p_0$ is not a peak point for $A(\overline{\Omega})$. Suppose it is. Then there is a $f \in A(\overline{\Omega})$, such that $f(p_0) = 1$ and $|f| < 1$ on $\overline{\Omega} \setminus \{p_0\}$. If we set $f_0(z) = f(z, 1)$, it follows that

$$(5) \qquad f_0(z_0) = 1, \quad |f_0(z)| < 1 \quad \forall z \in K \setminus \{z_0\}.$$

In view of Lemma 1 (ii), there is a nontrivial probability measure $\mu$ on $K$ such that

$$(6) \qquad |g(z_0)| \leq \int_K |g(z)| \, d\mu, \quad \forall g \in \mathcal{O}(K).$$

By Lemma 3, $f_0$ can be approximated by functions in $\mathcal{O}(K)$. Therefore, it follows from (6) that

$$(7) \qquad |f_0(z_0)| \leq \int_K |f_0(z)| \, d\mu < 1.$$

Since $\mu$ is not a point mass, (7) contradicts (5). This completes the proof of Theorem 1.

Now we prove that the same B-regular domain as in Theorem 1 also serves as an example in Theorem 2. First, we may assume that the point $X \ni \zeta_0 = 0$ by a simple Möbius transform (note that $\zeta_0$ is an interior point of the unit disc).

Let $X_t = \{\zeta \in X : |\zeta| \leq t\}$ for $t > 0$. Since $0 = \zeta_0$ is not a peak point for $\mathcal{R}(X)$, we claim that it is not a peak point for $\mathcal{R}(X_t)$ either, for any $t > 0$. Assume the contrary. Then there is a function $g \in \mathcal{R}(X_{t_0})$ for some $0 < t_0 < 1$ such that $g$ peaks at $0$. Choose $t_1, t_2$ such that $0 < t_2 < t_1 < t_0$. Set $E_1 = X_{t_2}$ and $E_0 = \overline{X \setminus X_{t_1}}$. The normality of the algebra $\mathcal{R}(X)$ yields a function $r \in \mathcal{R}(X)$ with $r(E_i) = \{i\}$ for $i = 0, 1$. Define a function $h$ by setting $h = rg$ on $X_{t_0}$ and $h = 0$ otherwise. Then we have that $h \in \mathcal{R}(X)$ and $h = g$ on $X_{t_1}$. In particular, $h$ is also a local peak function for $0$. Then by the local maximum modulus principle of Rossi [RO2], we know that $0$ is also a peak point for $\mathcal{R}(X)$, which is a contradiction. As a consequence, Bishop's theorem again implies that there is a nontrivial representing measure $\sigma_t$ for $0$ on $X_t$ for the algebra $\mathcal{R}(X_t)$, for any $t > 0$. By pushing forward the measure as before, we obtain a nontrivial representing measure $\mu_t$ for $z_0 = \Phi(0)$ on $K_t = \Phi(X_t)$ with respect to the algebra $\mathcal{O}(K_t)$. Now if $p_0 = (z_0, 1)$ is a local peak point for $\Omega$, then there is an open neighborhood $U$ of $p_0$ and a local peak function $f \in A(\Omega \cap U)$. Observe that $f$ can be approximated by functions $\{f_j\}$ that are holomorphic in some open neighborhood of $U \cap (K \times \{1\})$ by the same reasoning as before. It follows then that the restrictions of $f_j$'s to $K$ will be holomorphic in some neighborhood of $z_0$ in $K$ which contains $K_t$ for some sufficiently small $t > 0$. Finally, applying the same arguments as in the last part of the proof of Theorem 1 to $K_t$ and $\mu_t$, one can easily obtain that $|f(z, 1)| = 1$ for $\mu_t$-a.e. points in $K_t$, which is impossible because $f$ is a local peak function for $z_0 \in K_t$ and $\mu_t$ is not a point mass on $K_t$.

Department of Mathematics, Texas A&M University, College Station, TX 77843
*E-mail address*: jyu@math.tamu.edu